\documentclass[11pt]{amsart}
\usepackage{amsmath,amssymb}
\newtheorem{thm}{Theorem}

\newtheorem{lm}{Lemma}
\newtheorem{prop}{Proposition}

\newcommand{\be}{\begin{equation}}
\newcommand{\ee}{\end{equation}}
\newcommand{\bea}{\begin{eqnarray}}
\newcommand{\eea}{\end{eqnarray}}
\newcommand{\ba}{\begin{array}}
\newcommand{\ea}{\end{array}}
\newcommand{\beas}{\begin{eqnarray*}}
\newcommand{\eeas}{\end{eqnarray*}}

\def\Zhat{\widehat{Z}}
\def\e{\varepsilon}

\DeclareMathOperator{\Var}{Var}

\numberwithin{equation}{section}

\begin{document}

\title[Disordered exclusion process]
{Behavior dominated by slow particles
 in a disordered asymmetric exclusion process}

\author{Ilie Grigorescu}
\address{Ilie Grigorescu\\
Department of Mathematics\\
University of Miami\\
Coral Gables, FL 33124-4250 } 
\email{igrigore@math.miami.edu}
 
\author{Min Kang}

\address{ Min Kang\\
  Department of Mathematics\\
         North Carolina State University\\
         Raleigh, NC 27695}
\email{kang@math.ncsu.edu}

 \author{Timo Sepp\"al\"ainen}
\address{Timo Sepp\"al\"ainen\\
 Department of Mathematics\\
 University of Wisconsin--Madison\\ 
 Madison, WI 53706-1388}
\email{seppalai@math.wisc.edu}

%\date{\today}
\date{February, 2003}


%\shorttitle{Disordered exclusion process}

\begin{abstract}
We study the large space and time scale behavior of a 
totally asymmetric, nearest-neighbor
 exclusion process in one dimension  with random jump
rates attached to the particles. When slow particles are
sufficiently rare the system has a phase transition. At low
densities there are no equilibrium distributions, and on the 
hydrodynamic scale the initial profile is transported 
rigidly. We elaborate this situation 
further by finding the correct order of the correction from the 
hydrodynamic limit, together with 
 distributional bounds averaged over the disorder. 
We consider two settings,  a macroscopically constant
low density profile and  the outflow from a 
large jam. 
\end{abstract}

 \thanks{
Research of I.G. and M.K. partially supported by
the  Summer Internship Program in Probability
and Stochastic Processes, UW--Madison,  NSF grant DMS-0098605. 
Research of T.S.  partially supported by NSF grant DMS-0126775.}

\maketitle

\section{Introduction} 

We study a totally asymmetric, nearest-neighbor
 exclusion process
on the one-dimensional integer lattice $\bf{Z}$ with 
random rates attached to the particles. 
The process is studied through the labeled particle 
configuration. 
The particles are indexed  by integers in an 
increasing fashion. The position of particle $i$
at time $t$ is denoted by an integer-valued 
random variable  $\sigma_i(t)$. The exclusion 
rule stipulates that 
$
\sigma_i(t)< \sigma_{i+1}(t) 
$
for all $i\in\mathbf{Z}$ and all $t\ge 0$. 

At the outset each particle $\sigma_i$ receives its 
jump rate  $p_i$ which then remains fixed throughout the 
dynamics. The rates $\mathbf{p}=\{p_i\}$ 
are i.i.d.~random variables with common distribution $F$. 
 $F$ is supported on $(c,1]$
for some  $c>0$, and we take $c$ to be the left endpoint of the support
of $F$. In other words, $F(p)=0$ for $p<c$,   $F(p)>0$ 
for   $p>c$, and $F(1)=1$. We also assume $F(c)=0$ so no
particle has $c$ as its intrinsic jump rate. 

Once the rates have been fixed and an initial 
configuration $\sigma=(\sigma_i:i\in\mathbf{Z})$ specified,   the process
$\sigma(t)=(\sigma_i(t):i\in\mathbf{Z})$ evolves 
in the usual way: each particle $\sigma_i$ carries its own
Poisson clock of rate $p_i$, and whenever the clock rings, 
$\sigma_i$ advances one step to the right provided 
the next site to the right is vacant. 

It is also useful to consider the  gaps 
$\eta_i(t)=\sigma_{i+1}(t)-\sigma_i(t)-1$. The process
$\eta(t)=(\eta_i(t):i\in\mathbf{Z})$  is a zero-range process
with random rates attached to the spatial positions. 
The jump rule is that whenever a particle is present
at position $i$ [$\eta_i\ge 1$], one particle is moved
from $i$ to $i-1$ at rate $p_i$. We can also 
view this system as a series of tandem queues where 
queue $i$ is served at rate $p_i$, and customers departing
queue $i$ immediately join queue $i-1$. The gap variable 
$\eta_i(t)$ is the queue length 
  and the particle increment 
$\sigma_i(t)-\sigma_i(0)$ is the departure process
from queue $i$.

Fix the  rates $\mathbf{p}=\{p_i\}$. 
Given any $a\in[0,c]$, the product distribution $P^{\mathbf{p}}$ with 
geometric marginals 
\be
P^{\mathbf{p}}[\eta_i=k]=\Bigl(1-\frac{a}{p_i}\Bigr)
\Bigl(\frac{a}{p_i}\Bigr)^k\,,\quad k=0,1,2,\dotsc
\label{eq-dist}
\ee
is an invariant distribution for the gap process
$\eta(t)$.  
In this equilibrium, each particle motion is marginally
a Poisson process with rate $a$.
More precisely, for each $i$ the increment 
$\sigma_i(t)-\sigma_i(s)$ is Poisson with mean 
$a(t-s)$. 
 This is a consequence
of Burke's theorem from queueing theory, according to 
which the departure process of an M/M/1 queue in equilibrium
is a Poisson process. 

When $F$ is suitably chosen this model manifests a phase
transition. Here is a way to approach it. 
Given $a\in[0,c]$, the (annealed) mean gap in equilibrium is 
\[
u=\int E^{\mathbf{p}}[\eta_i]\,F^{\otimes\mathbf{Z}}(d\mathbf{p})
=\int_{(c,1]} \frac{a}{p-a}\,dF(p). 
\]
The common velocity $a$ of the particles cannot exceed $c$ 
because there are particles whose intrinsic rates 
come arbitrarily close to $c$. Thus the maximal mean 
gap $u^*$  is defined by letting  $a\nearrow c$, in other words
\be
u^*=\int_{(c,1]} \frac{c}{p-c}\,dF(p). 
\label{def-ustar}
\ee
If this integral is finite, there is a critical gap size
 $u^*<\infty$ such that the geometric product equilibrium  
distributions do not exist for mean gaps $u>u^*$. 
Equivalently, there is a positive critical density 
$\rho^*=(1+u^*)^{-1}$ for the exclusion particles such that the 
product equilibria for the gaps do not exist at low densities $\rho<\rho^*$. 
One interesting question is   the behavior of the system 
at low densities. 

 This system   attracted interest  in both
the theoretical physics and mathematics literature, starting 
from the mid--1990's. It appears that the invariant distributions
(\ref{eq-dist}) have been discovered several times independently. 
Among the early ones
was  Evans \cite{evan96, evan97} who derived the invariant
distributions for the disordered exclusion model  in both
continuous and  discrete time. 
Independently, Krug and Ferrari \cite{krug-ferr}
 studied the phase transition of the continuous-time model 
and  interpreted the results in various physical contexts
such as traffic flow and directed polymers. 
   In general, 
on the physics side there is  wide interest 
in particle systems as simple models of traffic flow 
and other ``single file'' systems. 
 We refer the reader to Nagel \cite{nage} for a review 
of particle systems in traffic modeling.
The state-of-the-art in 
 traffic modeling with exclusion type systems
 is the  Gray--Griffeath model
\cite{gray-grif},  
which is an exclusion
process whose jump rates depend on nearby sites. 

Returning to the disordered exclusion, 
on the mathematical side 
Benjamini et al.~\cite{benj-ferr-land} first 
proved hydrodynamic limits for several  
asymmetric exclusion and zero-range processes  with random rates. 
However, their assumptions specifically ruled out the phase
transition. 

A complete hydrodynamic limit theorem for the 
model studied here was proved by Sepp\"al\"ainen and Krug
\cite{sepp-krug}. For the case $\rho^*>0$
the result was the following. If the initial 
distributions have a  macroscopic profile   below  
 $\rho^*$, then on the hydrodynamic scale the initial 
macroscopic  profile
is rigidly translated at speed $c$.  In particular, if the 
system has initially a spatially homogeneous particle 
distribution   with density $\rho<\rho^*$ (such as ergodic
gaps with mean $u>u^*$),   a tagged
particle satisfies
\[
t^{-1}\sigma_i(t)\to c
\qquad\text{as}\ \ t\to \infty. 
\]

Subsequently Andjel et al.~\cite{andj-etal} 
proved a weak convergence result for the low density 
regime. Start the system so that the gaps are ergodic
with mean $u>u^*$. Then the gap process 
converges weakly to the maximal invariant 
distribution, in other words to the product distribution 
with marginals as  in (\ref{eq-dist}) with $a=c$. 

The hydrodynamic limit and the weak limit suggest the 
following picture. Let us follow particle $\sigma_0$ that
initially starts at the origin. The other particles are
distributed   so that the gaps are for example 
i.i.d.~with mean $u>u^*$, and then initially particle 
density is $\rho<\rho^*$.  As $t$ grows, particle $\sigma_0(t)$
experiences an increasing density around itself, 
and correspondingly its advance is slowed down. 
The reason is that $\sigma_0$ is part of an ever-growing
``platoon'' of particles, headed by an especially slow
particle. As this platoon catches up with  slower
platoons ahead of it, it grows and slows down even more.  
As $t\to\infty$, the particle density around $\sigma_0(t)$ 
approaches  the critical density $\rho^*$, and 
simultaneously his motion slows down to rate $c$. 
However, all this must happen at a scale below the 
hydrodynamic, because the hydrodynamic limit reveals only
the trivial final behavior. 

The purpose of this paper is to quantify the 
slowdown experienced by $\sigma_0(t)$ when the system 
starts at low density. Technically speaking we are 
seeking the next order term in the hydrodynamic limit. 
We find that by time $t$,   $\sigma_0(t)$ has traveled
a distance $ct+ w(t)t^{(\nu+1)/(\nu+2)}$ where $\nu>0$ is
an exponent characterizing the tail of $F(p)$ as $p\searrow c$, 
and $w(t)$ is a   random quantity, which 
becomes strictly positive and is tight as $t\to\infty$. 
We do not have 
a precise limiting distribution for $w(t)$.  
Our bounds suggest that for large $t$ the tail 
of $w(t)$ behaves like $\exp\bigl\{-C(u-u^*)^{-1} w^{2+\nu}\bigr\}$
for some constant $C$. These results are for annealed 
distributions, in other words for probabilities where the 
random rates have been averaged out. 

Following  the nonrigorous picture sketched above, 
proofs of the estimates proceed by bounding the rate of 
the slowest particle in a suitable range  ahead of $\sigma_0(t)$. 
The technical side of the proofs involves couplings of various 
kinds between several processes with different rates and/or 
initial distributions. 

We also address another question which is related, and partly
uses the same tools for the proof. When the exclusion process
starts with all sites in $(-\infty, 0]$ occupied and all
sites in $[1,\infty)$ vacant, on the hydrodynamic scale
there is a limiting macroscopic
profile supported on $(-\infty, c]$. 
In other words, the number $X_t$
 of particles that are in $(ct,\infty)$ at time $t$ must satisfy
$X_t=o(t)$. We find bounds on the true size of $X_t$. This 
question is not restricted to the situation where $u^*<\infty$.  
It makes sense whenever $F(c)=0$ because then every particle 
is attempting to jump at a rate strictly higher than $c$. 
Then presumably $X_t$ is unbounded as $t$ increases.

\section{The   results}

 The basic assumption is on the 
tail of $F(p)$ as $p\searrow c$. 
\be
\begin{split} 
&\text{There exist  constants $-1<\nu<\infty$  and $0<\kappa<\infty$
such that }\\
&\qquad\qquad\qquad
 \lim_{p\searrow c}\frac{F(p)}{ (p-c)^{\nu+1}} = \kappa. 
\end{split}
\label{ass1}
\ee
If the reader prefers   a concrete example, 
let $F$ have    density $f(p)=\kappa(\nu+1)(p-c)^\nu$
  on some interval $(c,c+\e)$. 
At $\nu=-1$ the distribution $F$ has a jump of size $\kappa$
at $c$, so there is  a positive density $\kappa$ of particles 
with minimal rate $c$. The behaviors we look at become 
simple. Values $\nu<-1$ are of course not possible.   
Recall the definition (\ref{def-ustar}) of the critical gap $u^*$. 
An integration by parts checks that, under
 assumption (\ref{ass1}),  $\nu>0$ is equivalent to 
$u^*<\infty$. 

First we look at the slowdown phenomenon in low density. 
We specify that particle 
$\sigma_0$ starts at the origin  [$\sigma_0=0$]. 
Initial  locations $(\sigma_i:i\ne 0)$ of the other particles 
are determined by taking the initial  gaps  $\{\eta_i\}$   i.i.d.~random
variables with common mean $u=E\eta_i>u^*$ and finite variance. 
Then set 
\[
\sigma_i=\sum_{j=0}^{i-1}\eta_j
\quad\text{for $i>0$, and}\quad 
\sigma_i=\sum_{j=i}^{-1}\eta_j
\quad\text{for $i<0$.}
\]
Our results are bounds on the ``annealed'' distributions of 
the quantities of interest. This means that while the 
process is run with fixed rates $\mathbf{p}=\{p_i\}$, we look 
at the average of all the processes for different choices
of $\mathbf{p}$, but with the  fixed initial distribution 
for $(\sigma_i)$. $P$ will denote this probability 
measure which represents the random choice
of rates, the random initial 
configuration $(\sigma_i)$, and the random exclusion evolution.

Notationally it is  convenient to use 
\[
\alpha=\frac1{\nu+2}\,,
\]
so that in particular the power  of the correction is 
\[
1-\alpha=\frac{\nu+1}{\nu+2}\,.
\]
Set also
\[
A(\nu)= \frac{(\nu + 2)^{\nu +2}}{(\nu +1)^{\nu+1}}.
\]

\begin{thm}\label{T} Assume {\rm (\ref{ass1})} with $\nu>0$. 
Let the initial gaps $\{\eta_i\}$ be  i.i.d.\ random
variables with common mean $u=E\eta_i>u^*$ and finite variance. 
The following bounds are valid for any $0<z<\infty$. 
\begin{equation}\label{ubc} \limsup_{t \to \infty}P\Bigl(
\frac{\sigma_{0}(t) - c t}{t^{1-\alpha}}
> z \Bigr)
\le \exp\Bigl\{ - A(\nu)^{-1}
\frac{\kappa  }{u-u^*} \, z^{\nu + 2}
  \Bigr\}
\end{equation}
and
\begin{equation}
\liminf_{t \to \infty} P\Bigl( \frac{\sigma_{0}(t) - c
t}{t^{1-\alpha}}
> z \Bigr) \ge \exp\Bigl\{ -  
 \frac{\kappa }{u-u^{*}}\, z^{\nu + 2}
      \Bigr\}\,.
\label{lbc}
\end{equation}
\end{thm}

Next we consider the situation where initially all sites
in $(-\infty, 0]$ are occupied by particles, and all sites
in $[1,\infty)$ are vacant. This could be thought of as an outflow
from a large jam. Now there is always a rightmost particle, 
so we label the particles with nonpositive integers
in increasing order. We drop the generic $\sigma$ notation, 
and for this special situation
 denote the  locations of the particles 
at time $t$ by 
\[
\dotsm <\xi_{-2}(t) < \xi_{-1}(t) <\xi_{0}(t).
\]
The initial locations are $\xi_{i}(0)=i$ for $i\le 0$. 
Particle $\xi_i$ jumps at rate $p_i$ independently
 drawn from distribution 
$F$.

This system has  a hydrodynamic limit which can be 
expressed  in terms of the empirical measure as follows:
for a compactly supported continuous test function $\phi$, 
\[
\lim_{t\to\infty} t^{-1}\sum_{i\le 0} \phi\bigl(t^{-1}\xi_i(t)\bigr)
=\int_{\mathbf{R}} \phi(x)r(x)\,dx
\]
almost surely. The limiting density $r(x)$ is supported on
$(-\infty, c]$. (The reader can find more information about 
the limit and $r(x)$ in \cite{sepp-krug}.) 
For the homogeneous exclusion with constant
rates 1 this is Rost's classical result \cite{rost}, with a piecewise
linear profile 
\[
r_1(x)=\begin{cases}
1, &x\le -1\\
\frac12(1-x), &-1<x\le 1\\
0, &x>1.
\end{cases}
\] 

The random rates produce the following qualitative difference 
with the homogeneous case. In the homogeneous case the 
lead particle $\xi_0(t)$ is a Poisson process of rate 1, and so its
location is $t+O\bigl({t}^{1/2}\bigr)$. In other words, its location 
coincides with the right edge of the hydrodynamic front. 
However, in the disordered system the lead particle 
is a Poisson process of rate $p_0$, which under assumption 
(\ref{ass1}) is strictly greater than $c$. Thus $\xi_0(t)$ 
and in fact  a large number of particles are
ahead of the hydrodynamic front whose right edge
at time $t$ is at $ct$. The second question we address
is to bound the number of these particles. 

Let $X_t$ be the number of particles that are beyond  
point $ct$ at time $t$, 
in other words
\[
X_t= 0\vee\,\sup \{ k\ge 1: \xi_{-k+1}(t)>ct\}.
\]

\begin{thm}\label{T2} Assume {\rm (\ref{ass1})} with $\nu>0$. 
 Then for all $b>0$, 
 \[
\limsup_{t\to\infty} 
P\{  X_t> bt^{1-\alpha}\}\le  
\exp\Bigl\{ - A(\nu)^{-1}  \kappa b^{\nu + 2}   \Bigr\}
\]
and 
 \[
\liminf_{t\to\infty} 
P\{  X_t> bt^{1-\alpha}\}\ge  
\exp\Bigl\{ -  A(\nu)   (1+u^*)^{\nu+1} \kappa b^{\nu + 2}   \Bigr\}
\]
\end{thm}

When $\nu\le 0$ we no longer have a finite critical 
gap size $u^*$. Theorem \ref{T} fails, not just because
$u>u^*$ is no longer possible, but because in equilibrium
$\sigma_0(t)$ is a Poisson process and has fluctuations
on the scale $t^{1/2}$. 

The phenomenon described by Theorem \ref{T2} is 
not restricted to $\nu>0$. With 
$-1<\nu\le 0$ it is still the case  that many particles
advance ahead of the hydrodynamic front, as no particle 
has the lower bound $c$ as its actual rate. 
  
For $\nu=0$ our result is the same as for $\nu>0$ but 
with a logarithmic weakening in the lower bound.
This seems an artifact of our proof, so it is not clear
whether this is  the true
state of affairs. 
Note that at $\nu=0$ we have $\alpha=1-\alpha=1/2$, 
matching with diffusive fluctuations. 

\begin{thm}\label{T3} Assume {\rm (\ref{ass1})} with $\nu=0$. 
Let $\e>0$ be arbitrarily small and 
$0<a<\infty$ arbitrarily large. If
$0<b<\infty$ is large enough, then   
for all large enough $t$, 
 \[
P\{ at^{1/2}(\log t)^{-1} \le X_t\le bt^{1/2}\}\ge 1-\e.
\]
\end{thm}

  $X_t$ 
changes behavior for  $\nu<0$, and 
is of smaller order  than $O(t^{1-\alpha})$.  Unfortunately 
 we do not have  matching upper and lower bounds. 
  As $\nu\searrow -1$ ($\alpha\nearrow 1$)
the ratio of the upper and lower bound exponents becomes one. 

\begin{thm}\label{T4} Assume {\rm (\ref{ass1})} with $-1<\nu<0$. 
Let $\e>0$.  If $0<b<\infty$ is large enough, then   
for all large enough $t$, 
 \[
P\{ X_t \le b t^{(1+\nu)/2}\} \ge 1-\e.
\]
If $0<a<\infty$ is small enough, then   
for all large enough $t$, 
 \[
P\{ X_t \ge a t^{(1+\nu)/(3+\nu)} \} \ge 1-\e.
 \]
\end{thm}

The upper bound is on the boundary of conflicting with 
Gaussian fluctuations of the Poisson clocks.
  For large $b$, with high probability the slowest particle 
among $bt^{(1+\nu)/2}$ particles  has rate at most 
$c+qt^{-1/2}$ for a small $q>0$. Consequently the 
number of jump attempts experienced by this slow particle 
by time $t$ is Poisson with
 mean $ct+qt^{1/2}$. This can be brought below $ct$ 
by a fluctuation of order $t^{1/2}$ in the clock.   Thus
there is some chance that this particle does not reach $ct$
by time $t$. To improve the probability to $1-\e$
we choose $b$ and $q$ so that there is a large enough 
number of slow particles.  The lower bound meets this ``Gaussian border''
only in the limit $\nu\searrow -1$.

\section{Variational representations} 

In this section we run through notions which have been 
elaborated elsewhere \cite{sepp-krug}. 
The purpose is to establish the conventions
followed in this paper which in some cases deviate slightly from 
those used before. 
Let an arbitrary initial configuration $\sigma=\{\sigma_i\}$
be given, random or deterministic. Fix the rates $\{p_i\}$. 
The process $\sigma(t)=\{\sigma_i(t)\}$ is  constructed 
with the usual graphical representation, by attaching a rate
$p_i$ homogeneous Poisson process $N_i=(N_i(t):t\ge 0)$
 to each particle $\sigma_i$. 

 Construct an auxiliary family $\{\zeta^i(t)\}$ of exclusion
processes by stipulating that at time $t=0$ their initial 
locations are
\[
\zeta^i_j(0)=\sigma_i+j\qquad\text{for $j\le 0$.}
\]
Only particle indices $j\le 0$ are used for the auxiliary 
processes. The jumps of the particles $\zeta^i_j$ are defined
by
\[
\text{$\zeta^i_j$ attempts to jump whenever Poisson clock 
$N_{i+j}$ rings.}
\]
This translation of the index of the clock has the effect 
that for any fixed $k$,
 particles $\{\sigma_k, \zeta^i_{k-i}: i\ge k\}$ make jump
attempts at the same times, namely when clock $N_{k}$ rings.

Process $\zeta^i(t)$ has initially all sites in $(-\infty, \sigma_i]$
 occupied and all sites in $[\sigma_i+1,\infty)$  vacant. 
From this observation one can see that the variational 
equation
\be
\sigma_k(t)=\inf_{i:i\ge k} \zeta^i_{k-i}(t)
\label{s-var-1}
\ee
is valid at $t=0$. Then one proves it by induction on jumps
for all times $t$. 

In Theorems \ref{T2}--\ref{T4} we consider the system 
$\xi(t)$ that starts exactly as $\zeta^i(t)$ but centered at 
the origin. Let 
\[
\xi^i_j(t)=\zeta^i_j(t)-\sigma_i. 
\]
Then the processes $\xi^i(t)$ are copies of $\xi(t)$, 
 except that the rates $\{p_i\}$ have 
been shifted in space. Of course this does not affect the distribution
of $\xi^i(t)$ when the rates are averaged out. We will find it convenient
to use the variational equality (\ref{s-var-1}) also in the form
\be
\sigma_k(t)=\inf_{i:i\ge k} \{\sigma_i+\xi^i_{k-i}(t)\}.
\label{s-var-2}
\ee

Exclusion processes can  be represented by interface 
processes. Suppose an interface process is given in 
terms of a height function $i\mapsto h_i(t)$
from $\mathbf{Z}$ into $\mathbf{Z}$. This means that at time $t$
the interface is the graph of the function $h(t)$, so that 
$h_i(t)$ is the vertical coordinate of the location of the 
interface over site $i$. We impose the condition 
$h_i\le h_{i+1}$ on admissible height functions. 
Dynamics are defined  by stipulating that 
if $N_i(t)=N_i(t-)+1$, then 
\[
\text{$h_i(t)=h_i(t-)+1$, provided $h_i(t-)\le h_{i+1}(t-)-1$.}
\]
In other words, height $h_i$ jumps up at rate $p_i$, provided
it does not go above   its right neighbor. 
Obviously, we can map between $\sigma(t)$ and $h(t)$ by 
\[
\sigma_i(t)=h_i(t)+i.
\]
Precisely speaking,  if the processes $\sigma(t)$ and $h(t)$
are coupled so that  this equality is true at $t=0$, then it 
remains true for all $t\ge 0$.

The gap process $\eta(t)=\{\eta_i(t)\}$ is defined in terms 
of these processes by
\[
\eta_i(t)=\sigma_i(t)-\sigma_{i-1}(t)-1=h_i(t)-h_{i-1}(t).
\]

The variational equation for the height process takes this 
form. Let 
$Z^i(t)$ be an interface process with these properties: 
initially 
\[
\text{$Z^i_j=0$ for $j\le i$, and $Z^i_j=\infty$ for $j> i$. }
\]
Dynamically, 
\[
\text{$Z^i_j$ takes its jump commands from Poisson process $N_j$,
for all $i$ and $j$. }
\]
Then 
\be
h_k(t)=\inf_{i:i\ge k}\{ h_i + Z^i_k(t)\}. 
\label{s-var-3}
\ee
There is no translation in (\ref{s-var-3}) because  each 
column
of the height processes $h(t)$ and $Z^i(t)$ reads the same clock. 
Since $\sigma_0(t)=h_0(t)$,   we can use the 
variational formula 
\be
\sigma_0(t)=\inf_{i:i\ge 0}\{ h_i + Z^i_0(t)\} 
\label{s-var-4}
\ee
 in the proof of Theorem \ref{T}
where we follow the evolution of $\sigma_0(t)$.

%As a final comment, we mention that the variational 
%representation of the exclusion process and other processes
%can be formulated as linear mappings in (max,plus)-algebra
%\cite{grig-kang}. 

\section{Proof of Theorem \ref{T}}

We begin with the key lemma that points the way
to controlling the behavior of the system by looking 
at the slowest rate in a suitable range of indices. 
 For fixed positive $q_{1}$ and $q_{2}$,
and a  positive real parameter $N$, 
let
\begin{equation}  J(N)= \inf\{i\ge 0:
p_i\le c + q_{2}N^{-\alpha} \} 
\label{def-J(N)}
\end{equation}
and define the event 
\begin{equation}\label{E:D}  D(N)=\Big\{p_i >
c + q_{2}N^{-\alpha} \ \text{for $0 \le
i \le [q_{1}N^{1-\alpha}]$} \Big\}
=\{ J(N)> q_{1}N^{1-\alpha}\} \,.
\end{equation}
 
\begin{lm}  
Assume  {\rm (\ref{ass1})} and recall that the rates $\{p_i\}$
are i.i.d.~with common distribution $F$. 
For  fixed  $q_1,q_2>0$, 
\begin{equation} \label{p4-rate}
\lim_{N\to \infty} P \big( D(N)\big) = \exp \bigl\{ 
- \kappa  q_{1} q_{2}^{\nu +1} \, \bigr\}\,.
\end{equation}
\label{DN-lm}
\end{lm}

\begin{proof}   
Let $\delta>0$.  For $p$ sufficiently close to
$c$,
\begin{equation}\label{fk}
(\kappa - \delta)(p-c)^{\nu+1} \le F(p) \le 
(\kappa + \delta)(p-c)^{\nu+1} \,.
\end{equation}
Due to the independence of the rates $p_{i}$, we have
\begin{equation} \label{p1-rate}
P \bigl(D(N)\bigr)\ =\
  \bigl(1-F(c + q_{2}N^{-\alpha}) 
\, \bigr)^{[q_{1}N^{1-\alpha}]} \,.
\end{equation}
This yields the upper and lower bounds
$$
\Bigl(\, 1 -  (\kappa \pm \delta)q_2^{\nu+1}N^{-\alpha(\nu+1)}\,
  \Bigr)^{[q_{1} N^{1 -\alpha }]}
$$
for $ P \bigl( D(N) \bigr)$. Let $N \to \infty$ and then $\delta \to
0$ to obtain the limit (\ref{p4-rate}).
\end{proof}

\subsection{Proof of the  upper bound in Theorem \ref{T}}

The upper bound (\ref{ubc})  follows from this proposition. 

\begin{prop} Suppose the initial gaps $\{\eta_i\}$ are an 
i.i.d.\ sequence with common mean $u>u^*$ and finite variance. 
Let $q_1, q_2>0$. Then for any $\delta>0$, 
$$
 \limsup_{t\to\infty} P[\sigma_0(t)\ge 
ct +(q_1(u-u^*)+q_2)t^{1-\alpha} +\delta t^{1-\alpha}]
\le
\lim_{t\to\infty} P(D(t)).
$$
\label{ub-prop}
\end{prop}

Before  proving  
Proposition \ref{ub-prop}, let us observe how it implies
the upper bound (\ref{ubc}). 
Together with Lemma \ref{DN-lm} the proposition gives
\[
 \limsup_{t\to\infty} P[\sigma_0(t)\ge 
ct + z t^{1-\alpha}]
\le\exp\bigl\{ -\kappa q_1q_2^{\nu+1}\bigr\}
\]
for any $q_1,q_2$ such that 
$
z=q_1(u-u^*)+q_2 +\delta.
$
Minimize the right-hand side of the inequality 
subject to this constraint on $q_1,q_2$. Then let $\delta\to 0$. 

The remainder of this section proves 
Proposition \ref{ub-prop}. 

\begin{lm} Consider an arbitrary process $\sigma(t)$. 
Let $K>1$. Then 
\[
\lim_{t\to\infty}P\Bigl[
\sigma_0(t)=\min_{0\le j\le Kt}\{ h_j+Z^j_0(t)\} \Bigr]=1.
\]
\label{Kt-lm}
\end{lm}

\begin{proof} From the definition of the process
$Z^{[Kt]}_i(\cdot)$, initially at time zero 
$Z^{[Kt]}_i(0)=0$ for $i\le [Kt]$. 
Variable $Z^{[Kt]}_{[Kt]}$ is the first to jump, after
which $Z^{[Kt]}_{[Kt]-1}$ may jump, then $Z^{[Kt]}_{[Kt]-2}$,
and so on. 
 Consequently the time $T$ 
when variable $Z^{[Kt]}_{0}$ takes its first jump up 
is a sum of independent exponential waiting times
with rates $p_{[Kt]}$, $p_{[Kt]-1}$, $p_{[Kt]-2}$, $\dotsc$, $p_0$. 
Let $\e>0$. 
Since each rate $p_i$ is bounded above by one, 
$T\le (K-\e)t$  with probability that vanishes 
exponentially fast as $t\to \infty$. If we take $0<\e<K-1$,
we conclude that 
\[
P\{Z^{[Kt]}_0(t)>0\}\to 0
\]
 exponentially fast.  

To prove the lemma, it suffices to   show that 
$Z^{[Kt]}_0(t)=0$ implies
\[
\sigma_0(t)=\min_{0\le j\le Kt}[ h_j+Z^j_0(t)].
\]
Let $i>[Kt]$. Then, since $Z^i_0\ge 0$ always and 
the height $h_i$ is nondecreasing in $i$,  
\begin{align*} 
h_i+Z^i_0(t) \ge h_i \ge h_{[Kt]} =h_{[Kt]}+Z^{[Kt]}_0(t).
\end{align*}
This shows that indices $i>[Kt]$ cannot contribute to the  
infimum in the variational formula. 
\end{proof}

\begin{lm} Consider two processes $\sigma$ and $\tilde\sigma$ 
whose initial gaps are i.i.d. with common mean 
$
E\eta_i=E\tilde\eta_i=u
$
and finite variances. Couple the initial configurations so 
that they are independent, but give the processes the same
rates $\{p_i\}$ and  the 
same Poisson clocks. Then  for $\delta>0$ 
$$
\lim_{t\to\infty} P[ \sigma_0(t)\ge \tilde\sigma_0(t)+\delta t^{1-\alpha}]
=0.
$$
\label{kolm-lm-1}
\end{lm} 

\begin{proof}
Let $K>1$ and define the event 
$$
A(K,t)=\left\{ \sigma_0(t)=\min_{0\le j\le Kt}[ h_j+Z^j_0(t)]  
\quad\mbox{and}\quad 
  \tilde\sigma_0(t)=\min_{0\le j\le Kt}[ \tilde{h}_j+Z^j_0(t)] 
\right\}.
$$
By Lemma \ref{Kt-lm}, 
 $P(A(K,t)) \to 1$ as $t\to \infty$. 
  On  $A(K,t)$, 
\beas
\sigma_0(t)&=&\min_{0\le j\le Kt}\{ h_j+Z^j_0(t)\} 
=\min_{0\le j\le Kt}\{ h_j -\tilde{h}_j +\tilde{h}_j +Z^j_0(t)\}\\
&\le& 
\tilde\sigma_0(t) +  \max_{0\le j\le Kt}\{ h_j -\tilde{h}_j\}.
\eeas
By Kolmogorov's inequality,
$$
P\left[  \max_{0\le j\le Kt}\{ h_j -\tilde{h}_j\} \ge 
\delta t^{1-\alpha}\right]
\le \frac{Kt \mbox{ Var}[\eta_1] + Kt \mbox{ Var}[\tilde\eta_1]}
{\delta^2 t^{2(1-\alpha)}}.
$$
As $2(1-\alpha)=2(\nu+1)/(\nu+2) >1$, this last expression
vanishes as $t\to\infty$. Consequently
\[
  P[ \sigma_0(t)\ge 
\tilde\sigma_0(t)+\delta t^{1-\alpha}]  
 \le P(A(K,t)^c) + 
P\left[  \max_{0\le j\le Kt}\{ h_j -\tilde{h}_j\} \ge 
\delta t^{1-\alpha}\right]. 
\]
gives the conclusion by letting $t\to\infty$. 
\end{proof}

Now define a particular mean $u$ initial system as follows. Fix 
a number $\bar{u}<u^*$, and let 
$\bar{a}$ be the equilibrium velocity corresponding to 
average gap $\bar{u}$, defined by
$$
\bar{u}=\int_{(c,1]} \frac{\bar{a}}{p-\bar{a}}\,dF(p).
$$
For each realization 
$\bf p$ of the rates, let $\{\bar\eta_i\}$ have the 
nonstationary geometric product equilibrium distribution 
$$
P^{\bf p} \left[\bar\eta_i=k\right]= \left(1-\frac{\bar{a}}{{p}_i}\right)
\left(\frac{\bar{a}}{{p}_i}\right)^k.
$$
Then 
$E\bar\eta_i=\bar{u}$, and the  $\{\bar\eta_i\}$ are i.i.d.\ when the 
random rates are averaged out. We chose $\bar{u}$ strictly less
than $u^*$ because then 
\[
E[{\bar\eta_i}^2]=\int_{(c,1]} 
\Bigl\{ 2\Bigl(\frac{\bar{a}}{p-\bar{a}}\Bigr)^2
+\frac{\bar{a}}{p-\bar{a}}\Bigr\}\,dF(p)<\infty.
\]
This finite variance is necessary so we can  apply the previous Lemma 
\ref{kolm-lm-1}. 
We cannot 
use  equilibrium gaps at mean $u^*$ because they
 have infinite variance if $0<\nu\le 1$. 

Let $\{\gamma_i\}$ be an i.i.d.\ sequence of nonnegative
integer valued random variables, independent of  $\{\bar\eta_i\}$,
and with common mean $E\gamma_i=u-\bar{u}$. Assume the  $\{\gamma_i\}$
have  finite variance. Define 
$$
\tilde\eta_i=\bar\eta_i+\gamma_i.
$$
Then  $\{\tilde\eta_i\}$ are i.i.d.\ with common mean $u$. 
Let $\tilde\sigma(t)$ denote the process with  
$ \tilde\sigma_0(0)=0$ and initial 
gaps  $\{\tilde\eta_i\}$. 

\begin{lm} For any $\delta>0$, 
$$
\limsup_{t\to\infty} P[\tilde\sigma_0(t)\ge 
ct +(q_1(u-\bar{u})+q_2)t^{1-\alpha} +\delta t^{1-\alpha}]
\le
\limsup_{t\to\infty} P(D(t)).
$$
\label{tildebar1}
\end{lm}

\begin{proof} Thinking of the zero-range process
 of the gap evolution, couple
the processes $\tilde\eta(t)=\{\tilde\eta_i(t)\}$ 
and $\bar\eta(t)= \{\bar\eta_i(t)\}$ via the basic
coupling, so that $\bar\eta_i(t) \le \tilde\eta_i(t)$ for all 
$i$ and $t$. This entails having 
$\tilde\sigma_i$ and $\bar\sigma_i$ read the same Poisson clocks
for each $i$. 

Let 
$$J(t)=\inf\{ i\ge 0: p_i \le c+q_2 t^{-\alpha}\}.
$$
$J(t)$ depends only on the rates. 
Since $\tilde\sigma_{J(t)}(t) -\tilde\sigma_{J(t)}(0)$ is 
stochastically dominated by a mean $ct+q_2 t^{1-\alpha}$
Poisson random variable, 
the event
$$
B_1(t)= \{ \tilde\sigma_{J(t)}(t) \le \tilde\sigma_{J(t)}(0)
+ct+q_2 t^{1-\alpha} +\delta t^{1-\alpha}/4\}
$$
satisfies $P(B_1(t)^c)\to 0$. 
Let
$$
B_2(t)= \{  \tilde\sigma_{J(t)}(0) \le J(t)(u+1)+ J(t)\delta/(4q_1) \}.
$$
By the weak law of large numbers, 
 $P(B_2(t)^c)\to 0$ because $J(t)\to\infty$ almost surely.
By the connection 
between particles $\tilde\sigma_i(t)$
and gaps $\tilde\eta_i(t)$, and by the  coupling with $\bar\eta_i(t)$, 
\beas
 \tilde\sigma_0(t)=  \tilde\sigma_{J(t)}(t) -
 \sum_{i=0}^{J(t)-1} \tilde\eta_i(t) -J(t)
\le  \tilde\sigma_{J(t)}(t) -
 \sum_{i=0}^{J(t)-1} \bar\eta_i(t) -J(t).
\eeas
By stationarity,   $\bar\eta(t)= \{\bar\eta_i(t)\}$ has the 
same distribution for all $t\ge 0$, under any fixed $\bf p$. 

Now combine the inequalities. On the event
$$
A(t)=\{\tilde\sigma_0(t)\ge 
ct +(q_1(u-\bar{u})+q_2)t^{1-\alpha} +\delta t^{1-\alpha}\},
$$
we have 
$$
\tilde\sigma_{J(t)}(t) \ge \sum_{i=0}^{J(t)-1} \bar\eta_i(t) +J(t)
+ ct +(q_1(u-\bar{u})+q_2)t^{1-\alpha} +\delta t^{1-\alpha}.
$$
Consequently on $A(t)\cap B_1(t)$ we have 
$$
\tilde\sigma_{J(t)}(0) \ge  \sum_{i=0}^{J(t)-1} \bar\eta_i(t) +J(t)
 + q_1(u-\bar{u})t^{1-\alpha} +\frac34\delta t^{1-\alpha}.
$$
Next,  on $A(t)\cap B_1(t)\cap B_2(t)$ we have 
$$
\sum_{i=0}^{J(t)-1} \bar\eta_i(t) \le J(t)u - 
 q_1(u-\bar{u})t^{1-\alpha} -\frac34\delta t^{1-\alpha} + 
\frac{J(t)\delta}{4q_1}.
$$
And finally, on the event $D(t)^c$, $J(t)\le q_1 t^{1-\alpha}$, 
and so as our last inequality, 
 on $A(t)\cap B_1(t)\cap B_2(t)\cap D(t)^c$ we have 
$$
\sum_{i=0}^{J(t)-1} \bar\eta_i(t) \le J(t)\bar{u} -\frac12\delta t^{1-\alpha}.
$$

To summarize, 
\beas
&& P[\tilde\sigma_0(t)\ge 
ct +(q_1(u-\bar{u})+q_2)t^{1-\alpha} +\delta t^{1-\alpha}]
 \le P(D(t)) + P(B_1(t)^c) \\
&&\qquad + P(B_2(t)^c) 
+ P\biggl( D(t)^c \cap\biggl\{ \,\sum_{i=0}^{J(t)-1} \bar\eta_i(t) \le 
J(t)\bar{u} -\delta t^{1-\alpha}/2 \biggr\}\biggr).
\eeas
The conclusion follows because on $D(t)^c$, $J(t)\le q_1 t^{1-\alpha}$
while still $J(t)\to\infty$, so the last probability vanishes
as $t\to\infty$. 
\end{proof}

Now we prove  Proposition \ref{ub-prop}. Fix $\bar{u}<u^*$ so that 
$$q_1(u^*-\bar{u})<\delta/4.
$$
 Define the processes 
$\tilde\sigma(t)$ and $\bar\sigma(t)$ as was done 
for Lemma \ref{tildebar1}. Couple all three processes
$(\sigma(t), \tilde\sigma(t), \bar\sigma(t))$ so that the 
initial gaps of $\sigma(t)$ are independent of the initial gaps 
of the other two, and all read the same Poisson clocks. 
By the choice of $\bar{u}$, 
\beas
&&P[\sigma_0(t)\ge 
ct +(q_1(u-u^*)+q_2)t^{1-\alpha} +\delta t^{1-\alpha}]\\
&&\qquad \le 
P[\tilde\sigma_0(t)\ge 
ct +(q_1(u-\bar{u})+q_2)t^{1-\alpha} +\delta t^{1-\alpha}/4]\\
&&\qquad\qquad + \; 
P[\sigma_0(t)\ge \tilde\sigma_0(t)  +\delta t^{1-\alpha}/4].
\eeas
Let $t\to\infty$ and apply the lemmas.

\subsection{Proof of the lower bound  in Theorem \ref{T}}

The lower bound will follow from proving this proposition.

\begin{prop}
Suppose the initial gaps are i.i.d.\ with common mean $u>u^*$ and 
finite variance. Given positive $q_1,q_2$ let
\[
r=\min\{ q_2,\,(u-u^{*})q_{1} \}.
\]
Then for any $\delta>0$, 
\begin{equation}\label{E:lbm}
\liminf_{t\to\infty}P\bigl\{ \sigma_{0}(t) \ge 
ct  +r t^{1-\alpha} -\delta t^{1-\alpha} \bigr\} 
\ge \exp\bigl\{- \kappa 
 q_{1}q_{2}^{\nu + 1} \bigr\}\,.
\end{equation}
\label{lb-prop}
\end{prop}

The lower bound (\ref{lbc}) will follow from this proposition 
the same way the upper bound (\ref{ubc}) followed from 
Proposition \ref{ub-prop}. Namely, for a given $z$ 
maximize the right hand side of (\ref{E:lbm})  
subject to $r-\delta=z$, and then let $\delta\to 0$.

To prove Proposition \ref{lb-prop}, we start with the 
variational equation and split it into two separate ranges.  
 $$ \sigma_{0}(t)= \inf_{j \ge
0}\{ h_{j}(0) + Z^{j}_{0}(t)\}= \min\{S_{1}(t),S_{2}(t)\}$$
where
\[
S_{1}(t)=\inf_{0 \le j \le q_{1}t^{1-\alpha}} \{ h_{j}(0) +
Z^{j}_{0}(t)\} \qquad\text{and}\qquad
 S_{2}(t)=\inf_{j >
q_{1}t^{1-\alpha}} \{ h_{j}(0) + Z^{j}_{0}(t)\}\,.
\]
We shall show that
\be
\label{lb-S1}
\liminf_{t\to\infty}P\bigl\{ S_1(t) \ge ct+q_2t^{1-\alpha}-
\delta t^{1-\alpha} \bigr\} \ge 
\lim_{t\to\infty}P\bigl(D(t))
 \ee
and 
\be
\lim_{t\to\infty} 
P\{ S_2(t)\ge ct +q_1(u-u^*)t^{1-\alpha}-\delta t^{1-\alpha} \}=1.
\label{lb-S2}
\ee
Together with Lemma \ref{DN-lm}, these imply (\ref{E:lbm}).

\subsubsection*{Proof of lower bound, part 1}

In this section we prove  (\ref{lb-S1}) for $S_1(t)$.

\begin{prop} Let  $q_1, q_2, \delta>0$.
 There exists an event $B(t)$ such that 
$P(B(t)^c)\to 0$  and 
$$
\{S_1(t) \ge ct +q_2t^{1-\alpha}-\delta t^{1-\alpha}\} 
\supseteq D(t)\cap B(t).  
$$
\label{S1-prop}
\end{prop}

Lower bound (\ref{lb-S1}) follows from this proposition. 
The rest of this subsection proves the proposition. 
Pick a further constant $q_3$ such that  
\[
0<q_3<q_2<q_3+\delta/4. 
\]
We shall couple $\sigma(t)$ with a faster process
$\hat\sigma(t)$ whose  jump rates $\hat{p}_i$ are given 
by
$$\hat{p}_i=p_i\vee(c+q_2N^{-\alpha}).$$
Process $\hat\sigma(t)$ will be in equilibrium so 
that each particle $\hat\sigma_i(t)$ jumps as a Poisson 
process with rate
$$
\hat{a}=c+q_3N^{-\alpha}.
$$
To achieve this, the gap process
 $\hat\eta(t)=\{\hat\eta_i(t)\}$ has to have the 
appropriate geometric product equilibrium distribution. 
Given $\bf p$, $\{\hat\eta_i\}$ are independent with 
geometric marginals
$$
P^{\bf p}\left[\hat\eta_i=k\right]= \left(1-\frac{\hat{a}}{\hat{p}_i}\right)
\left(\frac{\hat{a}}{\hat{p}_i}\right)^k\,, \ \ k=0,1,2,\dotsc 
$$
Note that this is sensible because $\hat{a}<\hat{p}_i$ for each $i$
by the assumption $q_3<q_2$. 
  The processes $\hat\eta(t)$ 
and $\hat\sigma(t)$ depend on
$N$, but we suppress this dependence from the notation. 

The mean gap for the $\hat\sigma(t)$ process is 
$$\hat{u}=E[\hat\eta_i]=\int_{(c,1]} \frac{\hat{a}}{\hat{p}-\hat{a}}\,dF(p).
$$ 

\begin{lm} The mean gap  $\hat{u}$ 
converges to $u^*$ as $N\to\infty$.
\label{uhat-lm}
\end{lm}

\begin{proof} The integral comes in two parts:
\begin{align*}
\hat{u}&=
\int_{(c,c+q_2N^{-\alpha}]} \frac{\hat{a}}{\hat{p}-\hat{a}}\,dF(p)
+
\int_{(c+q_2N^{-\alpha},1]} \frac{\hat{a}}{\hat{p}-\hat{a}}\,dF(p)\\
&=
\frac{c+q_3N^{-\alpha}}{(q_2-q_3)N^{-\alpha}} \cdot F(c+q_2N^{-\alpha})
+
\int_{(c,1]} \frac{c+q_3N^{-\alpha} }{p-c-q_3N^{-\alpha}}
\mathbf{1}_{(c+q_2N^{-\alpha},1]}(p)\,  \,dF(p).
\end{align*}
The first term on the last line vanishes as $N\to\infty$
by hypothesis 
(\ref{ass1}). To the second term we apply dominated convergence. 
The integrand converges to $c/(p-c)$ for each fixed 
$p\in(c,1]$, and satisfies the bound
\[
\frac{c+q_3N^{-\alpha} }{p-c-q_3N^{-\alpha}}
\mathbf{1}_{(c+q_2N^{-\alpha},1]}(p)
\le 
\frac{q_2}{q_2-q_3}\cdot\frac{c+q_3}{p-c}
\]
if $N\ge 1$. 
The last upper bound is integrable  under $dF(p)$, again by assumption
(\ref{ass1}). 
\end{proof}
 
For higher moments of $\hat\eta_i$ we develop
a bound. 

\begin{lm}  For $k\ge 1$
and $N\ge 4$,
$$E[\hat\eta_i^k]\le C N^{\alpha(k-1-\nu)^+}\log N.$$
  $C$ is a constant that depends on $k$ and all the 
other constants in the problem, but not on $N$. 
\label{etahatlm1} 
\end{lm}

\begin{proof} For a fixed  $\bf p$,  properties of a geometric
distribution give
$$
E^{\bf p}[\hat\eta_i^k] \le C_0 + C_1 (E^{\bf p}[\hat\eta_i])^k
\le C_0 + C_1  (\hat{p}_i-\hat{a})^{-k}
$$
for constants $C_0, C_1$ that depend on $k$. 
It remains to show 
$$
\int_{(c,1]} (\hat{p}-\hat{a})^{-k} \,dF(p)
\le  C N^{\alpha(k-1-\nu)^+}\log N.
$$
This integral is decomposed as 
\be 
(q_2-q_3)^{-k}N^{k\alpha} \int_{(c, c+q_2N^{-\alpha}]}  \,dF(p)
+
\int_{(c+q_2N^{-\alpha}, 1]}  ({p}-\hat{a})^{-k}  \,dF(p).
\ee 
Apply assumption (\ref{ass1}) to the first integral. 
  In the second
integral observe that 
$$
 ({p}-\hat{a})^{-k} \le q_4(p-c)^{-k} 
\quad\mbox{for}\quad 
q_4=\left(\frac{q_2}{q_2-q_3}\right)^k.
$$
Subsume the constants $q_i$    into 
constants $C_i$.  Thus the next upper bound is of the form 
\be 
C_0 N^{\alpha(k-1-\nu)} 
+ C_1 \int_{(c+q_2N^{-\alpha}, 1]}  (p-c)^{-k}\,dF(p).
\ee
Pick $C_2$ and $\delta>0$ so that $F(p)\le C_2(p-c)^{\nu+1}$ for
$c<p\le c+\delta$.
In the second term, the integral over $(c+\delta,1]$ is 
bounded by a 
constant. Over $(c+q_2N^{-\alpha}, c+\delta]$ 
integrate by parts. 
\begin{align*}
\int_{(c+q_2N^{-\alpha}, c+\delta]} \frac{dF(p)}{(p-c)^{k}}
&\le F(c+\delta)\delta^{-k}-
\int_{(c+q_2N^{-\alpha}, c+\delta]} F(p)\,d\bigl\{(p-c)^{-k}\bigr\}\\
&\le C_3+ C_2k 
\int_{c+q_2N^{-\alpha}}^{c+\delta} (p-c)^{\nu-k}\,dp.
\end{align*}
Consider  different cases for the   last integral.
 If $\nu>k-1$, it is bounded by a constant. If $\nu=k-1$, it is
bounded by
$ C_3+ C_4 \log N$. Finally in the case  $\nu<k-1$ it is 
bounded by $C_5  N^{(k-1-\nu)\alpha}$. In all cases the bound
given in the statement of the lemma works. 
\end{proof}

 Couple $\{\eta_i\}$ and 
$\{\hat\eta_i\}$ so that they are mutually independent. 

\begin{lm} For any $q>0$, $\delta>0$, 
$$
\lim_{N\to\infty} P\Bigl\{ \,\inf_{0\le j\le qN^{1-\alpha}}
[h_j-\hat{h}_j] <-\delta N^{1-\alpha} \Bigr\} =0.
$$
Note that the height function $\hat{h}$ changes with $N$
in the statement above. 
\label{kolm-lm-2}
\end{lm}

\begin{proof} Take $N$ large enough so that $u-\hat{u}>0$,
which can be achieved by Lemma \ref{uhat-lm} and the assumption
$u>u^*$. 
Then the probability in the statement of the Lemma is bounded
above by 
\begin{align*}
&P\Bigl\{ \inf_{0\le j\le qN^{1-\alpha}}
[h_j-\hat{h}_j-j(u-\hat{u})] <-\delta N^{1-\alpha} \Bigr\}\\
&\le \delta^{-2}N^{-2(1-\alpha)}\cdot q N^{1-\alpha}\bigl(
\Var[\eta_0]+\Var[\hat{\eta}_0] \bigr)
\end{align*}
where we used Kolmogorov's inequality.  By the previous
lemma $\Var[\hat{\eta}_0]\le CN^{\alpha(1-\nu)^+}\log N$ 
while 
$\Var[\eta_0]$ is a constant. As $\alpha(1-\nu)^+<1-\alpha$
for all $\nu>0$, the probability vanishes as $N\to\infty$. 
\end{proof}
 
Now we turn to $S_1(t)$. Consider first a fixed $t$.
Set
$N=t$, and as above construct the equilibrium process
$\hat\sigma(\cdot)$ with  rates $\hat{p}_i$. 
 Also, let $\Zhat^j$   denote the corner  
processes run with the $\hat{p}_i$ rates. 
 On the event $D(t)$,   we
have 
\beas
S_1(t)&=&\min_{0\le j\le q_1t^{1-\alpha}}\{ h_j+\Zhat^j_0(t)\} 
=\min_{0\le j\le q_1t^{1-\alpha}}\{ h_j -\hat{h}_j
+ \hat{h}_j +\Zhat^j_0(t)\} \\
&\ge& 
\min_{0\le j\le q_1t^{1-\alpha}}\{ h_j -\hat{h}_j\} +\hat\sigma_0(t)
\eeas
because 
$$
 \hat{h}_j +\Zhat^j_0(t) \ge \hat\sigma_0(t)
$$
for each $j\ge 0$. 
Consequently 
\beas
&& \left\{S_1(t) \ge ct +q_2t^{1-\alpha}-\delta t^{1-\alpha}\right\} \\
&\supseteq& D(t) \cap 
\biggl\{ \,\min_{0\le j\le q_1t^{1-\alpha}}\{ h_j -\hat{h}_j\} 
\ge -\tfrac12\delta t^{1-\alpha},\, \hat\sigma_0(t) \ge 
ct +q_2t^{1-\alpha}-\tfrac12\delta t^{1-\alpha}
  \biggr\} \\
&\equiv& D(t) \cap B(t),
\eeas
where the last identity means that the event $B(t)$
is defined by the previous expression in braces. 
For the complement
\begin{align*}
P\bigl(B(t)^c\bigr)&\le 
P\Bigl\{\, \min_{0\le j\le q_1t^{1-\alpha}}\{ h_j -\hat{h}_j\} 
< -\tfrac12 \delta t^{1-\alpha}\Bigr\}\\  
&\qquad\qquad + P \left\{ \hat\sigma_0(t) <
ct +q_2t^{1-\alpha}-\tfrac12\delta t^{1-\alpha} \right\}.
\end{align*} 
The   probabilities above  vanish as $t\to\infty$,
the first 
by  Lemma \ref{kolm-lm-2} above.
For the second probability, note that $\hat\sigma_0(t)$ is 
Poisson distributed with mean 
\[
\hat{a}t=ct+q_3t^{1-\alpha} > ct +q_2t^{1-\alpha}-
\tfrac14\delta t^{1-\alpha}.
\]
Since $1-\alpha>1/2$, the deviation $\tfrac14\delta t^{1-\alpha}$
has zero probability in the $t\to\infty$ limit. 
This completes the proof of Proposition \ref{S1-prop}.

\subsubsection*{Proof of lower bound, part 2}

In this section we complete the proof of Proposition 
\ref{lb-prop} by proving (\ref{lb-S2}). 

\begin{prop} Given $\e,\delta>0$, 
\[
P\{ S_2(t)\ge ct +(u-u^*-\delta)q_1t^{1-\alpha}\} \ge 1-\e
\]
for all large enough $t$.  
\end{prop}

\begin{proof} 
Let $\{\eta^*_i\}$ be the independent mean $u^*$ equilibrium gaps, so given 
$\bf p$, 
\[
P^{\bf p}\left[\eta^*_i=k\right]= \left(1-\frac{c}{{p}_i}\right)
\left(\frac{{c}}{{p}_i}\right)^k,\quad k\ge 0.
\]
Let $\sigma^*(t)$ be the equilibrium  process where particle  
$\sigma_0^*(t)$ is a rate $c$ Poisson process. Couple the 
processes $\sigma(t)$ and $\sigma^*(t)$ so that they read
the same Poisson clocks but their initial states are 
independent. 

By the Strong Law of Large Numbers, 
\[
M^{-1}(h_M-h^*_M)\xrightarrow[M\to\infty]{} u-u^*
\qquad\text{almost surely.}
\]
Note that here we do not need finite variance, which the 
$h^*$ height function would not possess if $0<\nu\le 1$. 
Shrink $\delta$ if necessary so that
 $0<\delta<u-u^*$. Pick $M_0=M_0(\delta,\e)$ such that  
\[
P\bigl\{ h_j-h^*_j \ge j(u-u^*-\delta/2)\;
\text{ for all $j\ge M_0$}
\bigr\} \ge 1-\e/2.
\]
Since  $1-\alpha>1/2$, there exists a $t_0$ such that
\[
P\{\sigma^*_0(t)\ge ct-q_1t^{1-\alpha}\delta/2\}\ge 1-\e/2
\]
for all $t\ge t_0$. 
Now with probability at least $1-\e$, for $t\ge t_0$ such that
$q_1t^{1-\alpha}>M_0$, 
\begin{align*}
S_2(t)&=\inf_{j\ge q_1t^{1-\alpha}}\{h_j+Z^j_0(t)\}
=\inf_{j\ge q_1t^{1-\alpha}}\{h_j-h^*_j+h^*_j+Z^j_0(t)\}\\
&\ge \inf_{j\ge q_1t^{1-\alpha}}\{h_j-h^*_j\} +\sigma^*_0(t)
\ge  q_1t^{1-\alpha}(u-u^*-\delta/2) 
        + ct-q_1t^{1-\alpha}\delta/2\\
&= ct +q_1t^{1-\alpha}(u-u^*-\delta).
\end{align*}
\end{proof}

\section{Proof of Theorem \ref{T2}}

We begin with  the upper bound. 
Let $b>0$, $0<\theta<1$, $0<\e<\theta b$ and $q_2=\theta b-\e$.
Let $\bar{p}$ be the minimal rate among 
$p_{-[bt^{1-\alpha}]},\dotsc,p_{-[\theta bt^{1-\alpha}]}$,
and $I$ an index such that $p_I=\bar{p}$. 
Let $Y(t)$ be a Poisson variable with mean 
$ct+q_2t^{1-\alpha}$. If 
$\bar{p}\le c+q_2t^{-\alpha}$,  $Y(t)$ dominates the number
of jump attempts particle $\xi_I$ experiences 
during time interval $[0,t]$. By the particle ordering, 
$\xi_{-[bt^{1-\alpha}]}(t)\ge ct$ implies
$\xi_I(t)\ge ct$, and thereby $\xi_I$ must have at least
$ct+\theta bt^{1-\alpha}$ jump attempts. We get the bound
\[
\begin{split}
P\{ X_t>bt^{1-\alpha}\} &\le P\{\xi_{-[bt^{1-\alpha}]}(t)\ge ct\}\\
&\le P\{ \bar{p}>c+q_2t^{-\alpha}\} + P\{ Y(t)\ge ct+\theta bt^{1-\alpha}\}.
\end{split} 
\]
Since $1-\alpha>1/2$, 
the last probability vanishes as $t\to\infty$. By Lemma \ref{DN-lm}
we get 
\[
\limsup_{t\to\infty} P\{ X_t>bt^{1-\alpha}\} \le 
\exp\bigl(-\kappa (1-\theta)b q_2^{\nu+1}\bigr).
\]
Let $\e\searrow 0$ so that $q_2\nearrow \theta b$, and then choose 
$\theta=(\nu+1)/(\nu+2)$. 

\hbox{}

The lower bound of Theorem \ref{T2} comes from the lower bound
of Theorem \ref{T}. Pick a density $u>u^*$, and let the initial
gaps $\{\eta_i\}$ be bounded i.i.d.\ random variables with mean $u$. 
 By the variational formula (\ref{s-var-2}), 
\[
\sigma_0(t)=\inf_{j\ge 0} \{ \sigma_j +\xi^j_{-j}(t)\}
\]
where $\xi^j(t)$ is a version of the $\xi(t)$ process 
with   translated rates. 
Let $b>0$, and then pick $\theta>b(u+1)$.  Let $j=[bt^{1-\alpha}]$. Then 
\[
\begin{split}
\xi^j_{-j}(t) \ge \sigma_0(t)- \sigma_j
\ge ct +\bigl( \sigma_0(t)-ct\bigr)   - \sigma_j. 
\end{split}
\]
The annealed
distribution of the process $\xi^j(t)$ is the 
same as that of $\xi(t)$. Consequently
\[
\begin{split} P\bigl\{  X_t>bt^{1-\alpha}\bigr\} \ge 
P\bigl\{\xi_{-[bt^{1-\alpha}]}(t)>ct\bigr\} \ge 
P\biggl\{ \frac{\sigma_0(t)-ct}{t^{1-\alpha}}> \theta \,,\,
\sigma_{[bt^{1-\alpha}]} < \theta t^{1-\alpha} 
\biggr\}. 
\end{split}
\]
By the law of large numbers  
$t^{-1+\alpha}\sigma_{[bt^{1-\alpha}]} \to b(u+1)$, and so
by Theorem \ref{T}, 
\[
\liminf_{t\to\infty} P\bigl\{  X_t>bt^{1-\alpha} \bigr\}
\ge 
\liminf_{t\to\infty} 
P\biggl\{ \frac{\sigma_0(t)-ct}{t^{1-\alpha}}> \theta\biggr\}
\ge
\exp\Bigl\{-\frac{\kappa}{u-u^*}\theta^{\nu+2}\Bigr\}.
\]
Maximize the last lower bound over $\theta$ and $u$ subject to 
$u>u^*$ and 
$\theta>b(u+1)$.

\section{Proof of Theorem \ref{T3}}

The argument for 
  the upper bound is similar to the previous one. Now 
  $\nu=0$   and   $1-\alpha=1/2$. 
Let $\e>0$ be small. 
By the central limit theorem we can
fix  a large $1<M<\infty$ so that, if
$Y(t)$ is   a Poisson  random variable   with mean $ct+t^{1-\alpha}$,
then 
\[
P[Y(t)\ge  ct+Mt^{1-\alpha}]\le \e/4
 \]
for all large enough $t$.  
Given $\e$ and $M$, choose $0<q_2<1<M <q<b$ 
so that 
\[
\exp(-\kappa q q_2^{\nu+1})\ge 1-\e/16 
\qquad\text{and}\qquad
\exp(-\kappa b q_2^{\nu+1})\le \e/16.  
\]
Let 
\[
J(t)=\inf\{i\ge 0: p_{-i}\le c+q_2t^{-\alpha}\}.  
\]
By  Lemma \ref{DN-lm} 
we have 
  $t_0<\infty$ so that 
\[
P\{qt^{1-\alpha}< J(t)< bt^{1-\alpha}\}\ge 1-\e/4
\]
for all $t\ge t_0$. Suppose this event happens. 
Then if $\xi_{-[bt^{1-\alpha}]}(t)\ge ct$, 
  also 
$\xi_{-J(t)}(t)\ge ct $, and  
particle  $\xi_{-J(t)}$ has had to cover  distance 
$ct+J(t)\ge ct+qt^{1-\alpha}$.  
The increment   $\xi_{-J(t)}(t)-\xi_{-J(t)}(0)$ is
stochastically bounded  
by   the variable $Y(t)$ defined above. 
 So for large enough $t$, 
\[
\begin{split}
P\{ X_t>bt^{1-\alpha}\}  &\le P\{\xi_{-[bt^{1-\alpha}]}(t)\ge ct\}
\le
P\{\xi_{-J(t)}(t)\ge ct\} +\e/4\\
&\le 
P\{Y(t)\ge  ct+qt^{1-\alpha}\}  +\e/4 \le \e/2. 
\end{split} 
\]

We prove the lower bound by comparison with a faster 
system in equilibrium. Let $0<a<\infty$
be fixed.  Given $\e>0$, pick $1<w<\infty$
large enough so that 
\be
P[Y(N)> EY(N)-wN^{1/2}]\ge 1-\e/4
\label{Y-ineq-1}
\ee
for large enough $N$, for a Poisson variable $Y(N)$ with
mean $cN+2wN^{1/2}$. Later we have to increase $w$ further. 

Let $q_2=4w$, and define faster rates by $\hat{p}_i=p_i\vee(c+q_2N^{-1/2})$.
Consider $N$ large enough to have $\hat{p}_i<1$. 
Let $\hat\sigma(t)$ be a process run with rates $\hat{p}_i$ and in 
equilibrium, so that $\hat{\sigma}_0(t)$ is a Poisson process
with rate 
\[
a=c+2wN^{-1/2}.
\]
 The gap process $\hat{\eta}(t)$ then 
has a product distribution 
with independent geometric marginals
\[
P^{\bf p}\left[\hat\eta_i=k\right]= \left(1-\frac{a}{\hat{p}_i}\right)
\left(\frac{a}{\hat{p}_i}\right)^k\,, \ \ k=0,1,2,\dotsc 
\]
The annealed mean gap is 
\[
u=E[\hat\eta_i]=\int_{(c,1]} \frac{a}{\hat{p}-a}\,dF(p)
\]
and the annealed variance is bounded as in 
\[
\Var[\hat\eta_i]\le E[\hat\eta_i^2]=
2\int_{(c,1]} \biggl(\frac{a}{\hat{p}-a}\biggr)^2\,dF(p)+u.
\]

\begin{lm} There is a  constant $C$ that depends only on the 
distribution $F$ such that for large enough $N$, 
\[
u\le C\log N  \qquad\text{and}\qquad \Var[\hat\eta_i]\le CN^{1/2}.
\]
\label{bd-lm-1}
\end{lm}

\begin{proof} First for the mean.  Integrate by parts, and use
assumption (\ref{ass1}) to pick $0<\delta<1$ such that
 $F(p)\le (\kappa+1)(p-c)$ for 
  $c<p<c+\delta$. Then note that 
\[
p-a\ge \frac{q_2-2w}{q_2}(p-c)=\tfrac12(p-c)
\qquad\text{for $p\ge  c+q_2N^{-1/2}$.}
\]
Carrying out these steps yields 
\[
\begin{split}
u&= \frac{aF(c+q_2N^{-1/2})}{c+q_2N^{-1/2}-a} +
    \int_{(c+q_2N^{-1/2},1]} \frac{a}{p-a}\,dF(p)\\
&=  \frac{aF(c+q_2N^{-1/2})}{c+q_2N^{-1/2}-a}\\
&\qquad\qquad\qquad +
\biggl\{ \frac{a F(1)}{1-a} - \frac{aF(c+q_2N^{-1/2})}{c+q_2N^{-1/2}-a}
 -  \int_{(c+q_2N^{-1/2},1]} F(p)  d\biggl(\frac{a}{p-a}\biggr)\,\biggr\}\\
&= \frac{a }{1-a} + a  \int_{c+q_2N^{-1/2}}^1 \frac{F(p)}{(p-a)^2}dp\\
&\le  \frac{a }{1-a} +4 (\kappa+1)a 
\int_{c+q_2N^{-1/2}}^{c+\delta} \frac{dp}{p-c} 
+4 a  \int_{c+\delta}^1 \frac{F(p)}{(p-c)^2}dp\\
&\le  \frac{a }{1-a} +4(\kappa+1) a 
\Bigl( \log\delta-\log q_2N^{-1/2}\Bigr) + 4a \delta^{-1}\\
&\le \frac{1+c}{1-c} + 2(\kappa+1)\log N +
 4 \delta^{-1}\\
&\le C\log N.
\end{split}
\]
In the second last step we took $N$ large enough so that 
\[
a=c+q_2N^{-1/2} \le \frac{1+c}2 \le 1. 
\]
If $N\ge 3$, in the last step we can take 
\[
C=\frac{1+c}{1-c}+ 2(\kappa+1)+  4 \delta^{-1}
\]
which depends only on the distribution $F$. 

Following the same pattern for $E[\hat{\eta}_i^2]$ shows that,
after  integration by parts,  
the main part is the integral 
\[
 a^2\int_{c+q_2N^{-1/2}}^1 \frac{F(p)}{(p-a)^3}\,dp
\le
8a^2(\kappa+1) 
\int_{c+q_2N^{-1/2}}^{c+\delta} \frac{dp}{(p-c)^2} 
+  8a^2   \int_{c+\delta}^1 \frac{F(p)}{(p-c)^3}dp.
\]
The desired bound follows as above. 
\end{proof} 

Let $\hat\xi(t)$ denote a $\xi$-type process run with rates $\hat{p}_i$. 
Let \[
j(N)=[aN^{1/2}(\log N)^{-1}].
\]
 From the variational  coupling (\ref{s-var-2})  we
have 
\[
\begin{split}
\hat{\xi}^{j(N)}_{-j(N)}(t) &\ge \hat{\sigma}_0(t)-\hat{\sigma}_{j(N)}\\
&= ct +2wN^{-1/2}t
+\bigl(\hat{\sigma}_0(t)-(ct +2wN^{-1/2}t)\bigr)   -\hat{\sigma}_{j(N)}.
\end{split}
\]
The processes  $\hat\xi(t)$ and  $\hat{\sigma}(t)$ depend on $N$
but we suppress this from the notation. Set time $t=N$. 
Note that when the random rates are averaged out, processes
$\hat{\xi}^{j(N)}(t)$ and  $\hat{\xi}(t)$  have the same 
distribution. 
We get this bound. 
\[
\begin{split} 
&P\bigl\{ \hat{\xi}_{-j(N)}(N)> cN \bigr\}\\
&\qquad\ge 
P\Bigl\{ \hat{\sigma}_0(N)-(cN +2wN^{1/2}) >-wN^{1/2}\,,\,
  \hat{\sigma}_{j(N)}<wN^{1/2}\Bigr\}\\
&\qquad\ge 
P\bigl\{ \hat{\sigma}_0(N)-(cN +2wN^{1/2}) >-wN^{1/2}\bigr\}
- P\bigl\{  \hat{\sigma}_{j(N)}\ge wN^{1/2}\bigr\}.
\end{split}
\]
The next to last probability is at least $1-\e/4$ for large $N$
by (\ref{Y-ineq-1}). 
It remains to show that the last probability vanishes as 
$N\to\infty$. 
From the annealed perspective $\hat{\sigma}_{j(N)}$ is a sum 
of i.i.d.'s, so its mean  and 
variance are bounded as follows: 
\[
E\hat{\sigma}_{j(N)}=j(N)(u+1) \le CaN^{1/2} 
  \quad\text{and}\quad 
\Var[\hat{\sigma}_{j(N)}]=j(N)\Var[\hat{\eta}_0]\le CaN(\log N)^{-1}.
\]
At this point we need to increase our original choice
of $w$ to guarantee that  
$w>2Ca$ where $a$ is given in the beginning of the proof and $C$
is the constant that appears in Lemma \ref{bd-lm-1}. Then 
Chebychev's inequality gives
\[
P\bigl\{  \hat{\sigma}_{j(N)}\ge wN^{1/2}\bigr\}
\le 
P\bigl\{  \hat{\sigma}_{j(N)}\ge E\hat{\sigma}_{j(N)}+ CaN^{1/2}\bigr\}
\le
\frac{ CaN(\log N)^{-1}}{ C^2a^2N}
\]
which vanishes as $N\to\infty$. We can conclude that for large
$N$, 
\[
P\bigl\{ \hat{\xi}_{-j(N)}(N)> cN \bigr\}>1-\e/3.
\]

Finally we make contact with $\xi(N)$. Given $q_2$ chosen above, 
pick $q_1>0$ small enough so that $\exp(-\kappa q_1q_2)>1-\e/7$. 
Let $D(N)$ be the event 
\[
D(N)=\bigl\{ p_i=\hat{p}_i \quad\text{for $-[q_1N^{1/2}]\le i\le 0$}\bigr\}. 
\]
By  Lemma \ref{DN-lm},  $P(D(N))>1-\e/6$ for large enough $N$. 
On the event $D(N)$, $\xi_i(t)=\hat{\xi}_i(t)$ for 
$-[q_1N^{1/2}]\le i\le 0$ and all $t\ge 0$, so in particular for $i=j(N)$ 
if $N$ is large enough. Consequently 
\[
\begin{split}
&P\{ X_N\ge  aN^{1/2}(\log N)^{-1}\} \ge 
P\{ \xi_{-J(N)}(N)>cN \} \\
&\qquad \ge P\bigl(\{ \xi_{-J(N)}(N)>cN\}\cap D(N) \bigr)
= P\bigl(\{\hat{\xi}_{-J(N)}(N)>cN\}\cap D(N) \bigr) \\
&\qquad \ge
P\bigl\{ \hat{\xi}_{-j(N)}(N)> cN \bigr\} - P(D(N)^c)
> 1-\e/3 -\e/6 =1-\e/2.
\end{split}
\]
This completes the proof of Theorem \ref{T3}.

\section{Proof of Theorem \ref{T4}}

\subsection{Proof of the upper bound of  Theorem \ref{T4}} 
The upper bound   is proved by comparison  
with independent particles. 
Let 
\[
\beta=\frac{1-\alpha}{2\alpha}=\frac{1+\nu}2.
\]
  For $b>0$ and  $q_2>0$, define
\[
K_t=\sum_{i=-[b t^{\beta}]+1}^0 \mathbf{1}\{p_i\le c+q_2t^{-1/2}\}.
\]

\begin{lm} Let $\{Y_j(t)\}$ be independent copies
of a Poisson random variable with mean
$ct+q_2t^{1/2}$, independent of the rates $\{p_i\}$ 
and thereby independent of $K_t$. Then given
$\e>0$, 
 if $q_2$ is small enough while
 $bq_2^{\nu+1}$ is  large enough, 
\[
P\bigl\{ \text{$Y_j(t)\ge ct $ for $1\le j\le K_t$}\bigr\} <\e
\]
for all large enough $t$. 
\label{Y-lm-T3}
\end{lm}

\begin{proof}
Fix a small $0<\delta<1/2$. Fix a positive integer $m$ large enough
so that $(\tfrac12+\delta)^m<\e/2$.
Pick $\e_0>0$  small enough
so that 
\[
P(\chi\ge -\e_0) <(1+\delta)/2 
\]
for a standard normal $\chi$. Let $q_2<\e_0\sqrt{c}$. 

 By assumption (\ref{ass1}),
for large $t$  $K_t$ is stochastically dominated by
 a binomial random variable with 
$[b t^{\beta}]$ trials and success probability 
$(\kappa+1)q_2^{\nu+1}t^{-(1+\nu)/2}$. Such a variable
converges weakly to a Poisson with mean 
$b(\kappa+1)q_2^{\nu+1}$ as $t\to\infty$. Thus we may
fix $b$ large enough   so that 
\[
P(K_t\le m)<\e/2
\]
for large enough $t$. 

%Since 
%\[
%\frac{q_2t^{1/2}}{\sqrt{ct+q_2t^{1/2}\,}} <\e_0, 
%\]
 By the choice of $q_2$ and the definition of $Y(t)=Y_1(t)$, 
\[
\begin{split}
P\bigl\{ Y(t)\ge ct\bigr\} 
\le 
P\biggl\{
\frac{Y(t)-EY(t)}{\sqrt{\Var Y(t)\,}}\ge -\e_0\biggr\}.
\end{split}
\]
Then by the central limit
theorem, for large enough $t$ 
\[
P\bigl\{ Y(t)\ge ct\bigr\} 
\le P(\chi\ge -\e_0) +\delta/2 <1/2+\delta. 
\]

Finally, as the $Y_j(t)$ are i.i.d.\ and 
independent of $K_t$, 
\[
\begin{split}
&P\bigl\{ \text{$Y_j(t)\ge ct $ for $1\le j\le K_t$}\bigr\} 
=
E\biggl[ \,\prod_{j=1}^{K_t} P\{Y_j(t)\ge ct\}\biggr] \\
&\qquad\le 
E\Bigl[ (\tfrac12+\delta)^{K_t}\Bigr] 
\le P(K_t\le m) + (\tfrac12+\delta)^m \le \e.
\end{split}
\]
\end{proof}

Fix $b$ and  $q_2$   so that the lemma is satisfied. 
Let 
\[
I_t=\{ -[b t^{\beta}]< i\le 0: p_i\le c+q_2 t^{-1/2}\}.
\] 
Once the rates $p_i$ have been chosen according to distribution
$F$ and $I_t$ determined,  give each index $i\in I_t$ 
an independent Poisson process $N_i(\cdot)$ 
of rate $c+q_2t^{-1/2}$. Thin $N_i(\cdot)$ appropriately
to get the correct rate $p_i$. These thinned processes
 are the Poisson clocks 
for indices $i\in I_t$. Meanwhile, give the other indices 
 their independent Poisson clocks. This way we can claim that
for each $i\in I_t$, the number
of jump attempts experienced by particle $\xi_i$ during 
$(0,t]$  is bounded above by 
the mean $ct+q_2t^{1/2}$ Poisson variable  $N_i(t)$ that is
independent of the rates $p_i$. 

Suppose
$
\xi_{-[b t^{\beta}]+1}(t)\ge ct. 
$
By the particle ordering, $\xi_i(t)\ge ct$ for all $i\in I_t$, 
which implies that $N_i(t)\ge ct$ for all $i\in I_t$. 
By the lemma above this event has probability less than 
$\e$ for large $t$. 
To summarize, we have shown 
that for an arbitrary $\e>0$, $b$ can be chosen so that 
\[
P\{ X_t\ge b t^{\beta}\} <\e
\]
for large enough $t$.

\subsection{Proof of the lower bound of Theorem \ref{T4}}
For an exclusion process
with constant rates $r$, for any $a>0$ and  $0<\gamma<1$, 
\be
\lim_{t\to\infty} \frac{\xi_{-[t^\gamma a]}(t)-rt}{(rt)^{(1+\gamma)/2}}
=-2\sqrt{a}\qquad\text{in probability.}
\label{gl-wh}
\ee
This statement is a consequence of a limit proved by Glynn and
Whitt \cite{glyn-whit} and the explicit computation of the value
$2$ on the right hand side first done in 
 \cite{sepp97queue}. See Lemma 4.1
in \cite{sepp97queue} for the derivation of (\ref{gl-wh}) from 
\cite{glyn-whit}. (But note that the process $\xi$ in 
\cite{sepp97queue} is not the same as $\xi$ in the present paper.)

Let 
$
\beta=({3+\nu})^{-1}. 
$
Let $0<a<\infty$ and $q=2\sqrt{a}+2$. Use assumption
(\ref{ass1}) exactly as in the proof of Lemma \ref{DN-lm}
to show that, given $\e>0$, if $a$ is small enough, then 
for large enough $t$  
\[
P\{ \text{$p_i\ge c+qt^{-\beta}$ for $-[at^{\beta(1+\nu)}]\le i\le 0$}\}
\ge 1-\e/2.
\]
On this event $\xi_{-[at^{\beta(1+\nu)}]}(t)$ is bounded below 
by $\tilde{\xi}_{-[at^{\beta(1+\nu)}]}(t)$ where $\tilde\xi(t)$ is a 
process whose clocks ring at  constant rate $c+qt^{-\beta}$. 
For $\tilde\xi(t)$ 
 (\ref{gl-wh}) gives the following bound: 
 for large $t$   with 
probability at least $1-\e/2$, 
\[
\tilde{\xi}_{-[at^{\beta(1+\nu)}]}(t) \ge ct+qt^{1-\beta}
-2\sqrt{a}(ct+qt^{1-\beta})^{(1+\beta(1+\nu))/2}
- t^{(1+\beta(1+\nu))/2}
>ct.
\]
The last lower bound by $ct$  followed from $1-\beta 
=(1+\beta(1+\nu))/2$ and the choice of $q$. 

We have shown that, given $\e>0$ and 
a small enough $a>0$,  
then for large enough $t$, the inequality  
${\xi}_{-[at^{\beta(1+\nu)}]}(t) >ct$ holds
with probability at least $1-\e$.  This inequality   implies
 $X_t\ge at^{\beta(1+\nu)}$.  

\bibliographystyle{abbrv}
\bibliography{refs}

\end{document}